# Time Series Clustering
## Consistency of DBSCAN for clustering multivariate time series.

Nicholas Waltz, University of Cambridge

**Table of Contents**

# Abstract

Economic policy and research rely on the correct evaluation of the billions of high-frequency data points that we collect every day. Consistent clustering algorithms, like DBSCAN, allow us to make sense of the data in a useful way. However, while there is a large literature on the consistency of various clustering algorithms for high-dimensional static clustering, the literature on multivariate time series clustering still largely relies on heuristics or restrictive assumptions. The aim of this paper is to prove a notion of consistency of DBSCAN for the task of clustering multivariate time series.

# 1. Introduction and Motivation

The applications of time series clustering are vast. They range from the evaluation of weather data (e.g., Li *et al.*, 2016) to clustering medical data for disease classification (e.g., Malik *et al.*, 2019). Clustering makes it possible to sift through the enormous amounts of data that we collect every (nano-)second. Especially the fields of Finance and Economics rely on the correct evaluation of billions of data points. But which ones belong together? Clustering can achieve tasks like grouping a dataset on 1,000,000 stocks into meaningful portfolios or the classification of billions of income and consumer data.

The need for consistent time series clustering subroutines in Economic research is often overlooked, in fact, barely any empirical papers include clustering. However, when dataset selection and grouping is done by a researcher, this can and does result in bias and spurious regression phenomena, either from the omission or false inclusion of data points, be it out of naivety or due to data snooping. Consistent clustering algorithms are not only able to alleviate such bias, they are able to detect patterns in datasets that a human researcher would never be able to, especially given the billions of datapoints that are recorded every day.

In order for a clustering algorithm for time series to be useful, it needs to yield a consistent estimate of the clustering structure in the (hypothetical) population. However, while there is a large literature demonstrating the consistency of various non-parametric clustering algorithms in the static case (e.g., Hartigan, 1977), this problem has proved more difficult in the time series case. In fact, the majority of the time series clustering literature lacks such statistical rigour and instead relies on heuristics or finite sample performance measures, such as simulation studies.

Where there are asymptotic results, authors usually impose restrictive modelling assumptions (e.g., Jacques & Preda, 2014). This paper will aim to show that under relatively mild assumptions, we can apply results from the static case to the time series case.

As opposed to the static problem, clustering time series non-parametrically comes with additional difficulties. For one, it is unclear what we are comparing. In evaluating time series, we require a structural assumption, namely, whether we are dealing with functional data (functional data analysis, see Wang *et al.* (2015) for a review), or stochastic processes. This paper will make the functional data assumption, as most Economic indicators exhibit functional characteristics, e.g., business cycles, trade volumes, etc.

Secondly, given a structural assumption, it is unclear what we *should be* comparing. One approach, called raw time series clustering, involves treating the time series as vectors and comparing them directly as in the static case. However, given the same series twice, but out of phase, standard similarity measures, such as the Euclidian distance, may fail to accurately cluster the data. For this reason, authors usually advise against raw clustering (e.g., Jacques and Preda, 2013). In addition, the raw problem is rarely clearly defined. Reasons for this include incomplete data, warping, in-sample noise, and problem of choosing the subseries that should be compared.

For these reasons, time series require pre-treatment called filtering and registration. To deal with incomplete data and noise, kernel smoothers, such as the Nadaraya-Watson (Nadaraya, 1964; Watson, 1964) can be used, while motif-detection algorithms (see Fuchs *et al.*, 2009) can be used as a subroutine to deal with phase problems and subseries selection.

Given appropriate pre-treatment of the data, there are a lot of clustering algorithms that may be applied. By far the most popular clustering algorithm is k-means, or k-medoids, which has been studied extensively for the clustering of time series using a variety of similarity measures (e.g., Ashkartizabi & Aminghafari, 2017; Zambom *et al.*, 2018; Schmutz *et al.*, 2020). While intuitive and computationally simple, k-means has major shortcomings. For one, we are required to specify the number of clusters in advance, which in most practical cases, is an unknown quantity. And even given an estimator of the number of groups, there is no way to show that k-means will yield consistent clusters as computing the minimum of the k-means risk function is *NP*-hard (Dasgupta & Freund, 2009).

Other algorithms, most notably hierarchical clustering, k-nearest neighbours, self-organising maps and DBSCAN (Ester *et al.*, 1996) have better statistical properties and the literature on their clustering consistency in static problems, especially in conjunction with kernel density estimation, is well established.

The problem considered by this paper specifically, is the non-parametric clustering of multivariate functional data using DBSCAN, where our assumptions are suitable to an Economic setting. An example could be a data set containing company profits, employee numbers, and capital stock over time for a large number of companies. We abstract from registration and assume the data points are adequately pre-treated.

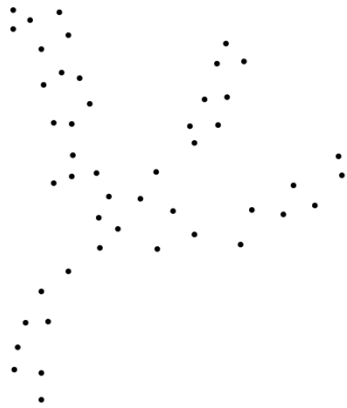

*Figure 1: How to cluster these 5 time series?*

The aim of this paper is to establish the clustering consistency in such a case. This paper will contain three contributions.
1. We will demonstrate that the multivariate problem can be flattened, and thus show that the usual concentration bounds, established in the literature on kernel density estimation, hold.
2. We will establish consistency, when the functional form first needs to be estimated from noisy and incomplete data.
3. We will include some discussion on dimensionality issues.

The rest of the paper is structured as follows. Section 2 provides a brief review of the literature on multivariate time series and clustering consistency. In Section 3, the problem and approach of this paper are specified and justified, after which the setup is provided in more formal language in Section 4. In Section 5, the main assertions of this paper are made and the consistency results are shown, and after a discussion of the practicalities of the problem in Section 6, we conclude the paper. For readability purposes, all proofs of Theorems, Lemmata, Corollaries and Remarks are provided at the end of the paper.

## 2. Review of relevant literature

### 2.1 Multivariate functional data clustering

The majority of the literature on time series clustering deals with the univariate case (Aghabozorgi *et al.*, 2015). As such, the clustering problem for univariate time series is well studied and different algorithms have been investigated. Particular attention is given to $k$-means. Because of this and due to curse of dimensionality problems, the majority of the multivariate functional data clustering literature employs some form of dimensionality reduction. One popular approach is model-based clustering, such as in Li *et al.* (2016), or Jacques and Preda (2014), where the clustering is done in parameter estimate space. As the parameter estimates are usually well-behaved, consistency is mostly immediate. The problem with such parametric approaches, however, is exactly that they are parametric. Not only does this mean that the clustering results are meaningless in the case of misspecification, a researcher can potentially assume the model that yields the most desirable result. For this reason, we should ideally employ a non-parametric approach. One standardised tool for non-parametric dimensionality reduction is functional principal component analysis (FPCA; e.g., in Ben-Hur & Guyon, 2003; Shang, 2013; Park & Ahn, 2017) which was first proposed for use in functional data analysis in Ramsay (1991). The idea of FPCA is to represent the functions in their eigenbasis and then conduct the clustering on the dominant modes of variation. The details of FPCA will not be discussed further in this paper but the reader is referred to Shang (2013) for a review. While using FPCA as a subroutine for dimensionality reduction has been shown to perform well in numerical studies (see Berrendero *et al*., 2011, Park & Ahn, 2017), it is more of a rule-of-thumb approach and the literature lacks consistency-type results. The main reason for this, is that in performing dimensionality reduction using FPCA, we inevitably lose information, albeit in the best way, as the eigenbasis expansion in FPCA accounts for more variation than any other basis expansion (Jones & Rice, 1992). However, it appears that in most practical applications, this bias does not matter. In high dimensional problems, where the clustering bias resulting from the loss of information, i.e., the mass of misclustered points, is comparatively small, misclustering rates will be negligible (Rinaldo & Wasserman, 2011).

Indeed, most clustering applications that use FPCA will be biased due to their small sample anyway. For consistency, however, we will need the full data.

## 2.2  Consistency results for static-data density clustering

The density clustering literature assumes an underlying data-generating density, such that the clusters lie in the denser regions of space. Much of the early consistency results for density clustering deal with the consistency of hierarchical clustering under different linkage algorithms. The earliest work on this can be found in Hartigan (1977), where a proof of the consistency of the single-linkage algorithm for estimating the cluster tree for densities on closed intervals in $\mathbb{R}$ is provided. Many of the concepts used in the density clustering literature are established in Hartigan (1981). Hartigan (1981) also demonstrates that the single linkage algorithm does not yield consistent estimates of the cluster tree in higher dimensions. Hierarchical clustering reconstructs the entire cluster tree for every level set. Subsequently, a lot of work has focussed on the theory of consistently estimating the clusters in a given level set (see Polonik, 1995; Tsybakov, 1997; Rigollet & Vert, 2009; Singh *et al.*, 2009; Rinaldo & Wasserman, 2011), often abstracting from the computational difficulties involved in detecting the clusters. Other work has focussed on improving on the shortcomings of the single-linkage algorithm with a particular emphasis on establishing convergence rates for the $k$-nearest neighbours approach (Maier *et al.*, 2009; Chaudhuri & Dasgupta, 2014). A particular area of focus in the clustering literature is the consistent estimation of the density using kernel density estimation (KDE). Given a consistent estimate of the density, clusters can be computed as the connected components of a given level set. This is also useful as density-based algorithms like hierarchical clustering, DBSCAN and k-nearest neighbours can all be related to KDE in some form. The theoretical basis for consistency in KDE can be found in Talagrand (1995) and Talagrand (1996), and later Giné and Guillou (2002), where uniform consistency and concentration bounds over all kernels with finite VC dimension are established. Later contributions like Kim *et al.* (2019) have further sharpened these rates. The contributions on KDE consistency have enabled proofs of the consistency of DBSCAN, mostly under the assumption of a Hölder-continuous data generating density, most notably in Wang *et al*. (2019) and Sriperumbudur and Steinwart (2012).

## 3. Preliminaries

### 3.1 Notation

In the following, let $\mu$ denote the Lebesgue measure. $\mathcal{B}(S)$ is the Borel $\sigma$-algebra on the topological space $(S, \mathcal{T})$, where for our purposes $\mathcal{T}$ is always the standard (Euclidian/$L^p$) induced topology. $B_d(x, r)$ is the Euclidian closed $d$-ball, where the index will be omitted, when the dimensionality is obvious from the context. Further, denote by $v_d = \mu(B_d(0,1)) = \frac{\pi^{d/2}}{\Gamma(\frac{d}{2}+1)}$ the volume of the Euclidian $d$-ball. $\|x\|$ and $\|x\|_\infty$ will denote the Euclidian and uniform norms, respectively. The uniform norm will be used both for functions and vectors (vectors can be interpreted as functions on finite subsets of $\mathbb{N}$, so there is no redundancy here). For two sets $A, B \subset S \subset \mathbb{R}^d$, $d(A, B) = \inf_{a \in A, b \in B} \|a - b\|$. $\bar{A}$ is the closure, $\dot{A}$ the interior of $A$ and $\partial A = \bar{A} \backslash \dot{A}$, the boundary of $A$. By $A^c = S \backslash A$, we denote the conjugate of $A$.

### 3.2 Density clustering and DBSCAN

We study clusters in the sense of Hartigan (1981). For data distributed according to a density $p(x)$ with support $S$, the $\lambda$-clusters $K_1^\lambda, K_2^\lambda, \ldots$ are the connected components of the level set $L(\lambda) \equiv \{x \in S : p(x) \geq \lambda\}$. In particular, $L(\lambda) = \cup_i K_i^\lambda$, and $K_i^\lambda \cap K_j^\lambda = \emptyset$, when $i \neq j$. Also, $L(0) = S$. Throughout we maintain that a density exists. For a given density $p$, the cluster tree $T$ is then the set of clusters for any given $\lambda$,

$$T(\lambda) = \{K_1^\lambda, K_2^\lambda, \ldots\}$$

The clusters in $T$ form a hierarchy in the sense that clusters get smaller as we increase $\lambda$ (Chaudhuri & Dasgupta, 2014). In particular, for $\lambda_1 \leq \lambda_2$,

1. $\forall K \in T(\lambda_2)$, there exists $K' \in T(\lambda_1)$, such that $K \subseteq K'$, and
2. For any $K \in T(\lambda_2)$ and $K' \in T(\lambda_1)$, either $K \subseteq K'$ or $K \cap K' = \emptyset$.

This means that the clusters are defined by the modes of $p$. Throughout, $\lambda$ is assumed to be a given quantity, and the aim is to show clustering consistency for arbitrary values of $\lambda$, which corresponds to a consistent estimate of the cluster tree $T$. There exist discussions on the choice of an optimal $\lambda$, given a sample, e.g., in Sriperumbudur and Steinwart (2012). We abstract from this. Consistency is shown in the sense of Hartigan (1981), and Chaudhuri and Dasgupta (2014).

**Definition 3.2.1.** Hartigan consistency *(Chaudhuri and Dasgupta, 2014). Say we construct an estimate $\hat{T}_n$ of $T$, given a sample $D_n = \{x_1, x_2, \ldots, x_n\} \subset S$. Let $A, A' \in L(\lambda)$ with $A \cap A' = \emptyset$ be disjoint sets in some level set, and $A_n, A'_n$ be the smallest clusters in $\hat{T}_n$ containing $A \cap D_n$ and $A' \cap D_n$ respectively. Then $\hat{T}_n$ is* Hartigan consistent *if whenever $A, A'$ lie in different elements of $T(\lambda)$, $\mathbb{P}(A_n \cap A'_n = \emptyset) \to 1$ as $n \to \infty$.*

One way to achieve a Hartigan-consistent estimate of $\hat{T}_n$ is to have a uniformly consistent estimator $\hat{p}_n$ of the density $p$ with $\|\hat{p}_n - p\|_\infty \to 0$. Then, we can let $\hat{L}_n(\lambda) = \{\hat{p}_n \geq \lambda\}$ and compute the connected components $\hat{T}_n(\lambda)$. Approaches for estimating the density in a uniformly consistent way are well studied in the clustering literature, specifically different kernel density estimators (among others Giné & Guillou, 2001; Rinaldo & Wasserman, 2011), but also histogram estimates (Steinwart, 2011). The problem is that finding the connected components is computationally difficult. To show that two points are in disjoint clusters, we have to show that every path between them fails to connect them.

DBSCAN gets around these computational difficulties, by considering balls around $D_n \cap \hat{L}_n$ instead of computing $\hat{L}_n$. This way the problem is agnostic about the underlying density, which is more tractable computationally and has been shown to perform very well for large samples (Ester *et al.*, 1996). Proofs of the consistency of DBSCAN for the static problem under different assumptions can be found in Rinaldo and Wasserman (2011), Sriperumbudur and Steinwart (2012), and Wang *et al.* (2019).

## 3.3 Intuition

We now move on to the main part of the dissertation. To motivate the density-based approach that follows, we consider the underlying structure of our clustering problem. In non-parametric functional data clustering, a common assumption on the data is that the observed data are of the form

$$y_i(t) = f_i(t) + \varepsilon_{it}$$

With suitable assumptions on the functions $f_i$ and (measurement) error processes $\varepsilon_{it}$, usually square-integrability and continuity for the functions and some form of independence for the shocks (Jacques & Preda, 2014). The usual procedure is then to reconstruct the functions $f_i$, from the noisy and incomplete data, mostly using some form of smoothing. The clustering task then consists of grouping the functions $f_i$ according to some similarity or distance measure, such as $L^p$ norms or radial basis functions.

In a sense, our data set is then a sample of the entire population of possible functions that belong to the category of data we observe. Since we expect the data to fall into groups of similar (close) functions, we expect certain time paths to be more likely than others. This means that we essentially expect a distribution over all time paths (see Figure 2). Therefore, under suitable framing, we can use the clusters-as-modes interpretation in density clustering: clusters are the dense regions of the time paths. Defining a meaningful probability measure over functions will prove a difficult task, however. Unless we impose restrictive assumptions, namely some form of countability or finiteness (e.g., finite VC dimension), most functional spaces are not measurable.

Upon closer inspection, this is a non-problem. While the underlying processes may be uncountably infinite dimensional in time space, computers will only compute similarity measures based on a countable, finite number of discrete point values (or estimates). Additionally, depending on the data, it may suffice to evaluate all the functions at a specific selection of time points to achieve correct clustering of the data. Then, the clustering problem essentially becomes one in $\mathbb{R}^n$, and we are able to apply consistency results from the literature.

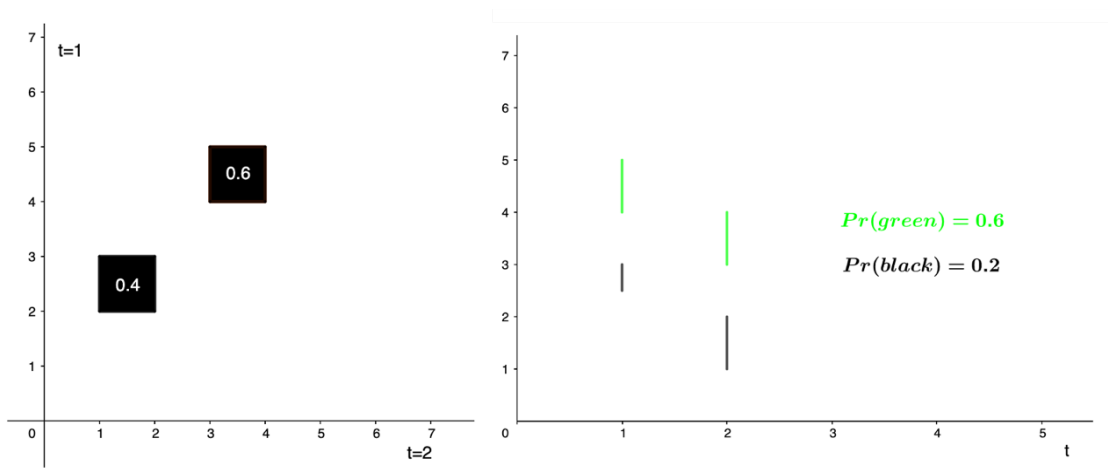

*Figure 2: Simple distribution and two examples of time path sets of univariate functions and their probability*

Ultimately, we can only cluster the data we are given, and the aim of this paper is therefore to prove consistency of DBSCAN given finite time samples. As this is how clusters are computed in practical applications (by digital computers), a consistency result for this approach is actually more powerful than assuming knowledge of the entire function, in the sense that it is more representative of reality.

## 4. Setup

In the following we consider the problem of clustering a set $D_n$ of $n$ multivariate time series $D_n = \{X_1, X_2, \ldots, X_n\}$, generated by a set of corresponding $s$-dimensional functions $f_1, f_2, \ldots, f_n$, whose components are defined on some interval, say $[0,1]$, are right-continuous, with a finite number of discontinuities, are locally $l_f$-Lipschitz on every continuous segment, and have an interval range, say $[0,1]$. Then we have

$$f_i(t) = \left(f_i^1(t), f_i^2(t), \ldots, f_i^s(t)\right)^T \in [0,1]^s \equiv \mathcal{S}$$

where $\forall t \in [0,1]$, $\lim_{c \to t^+} f_i^k(c) = f_i^k(t)$ and $f_i^k(t) \in [0,1]$. For clustering, we evaluate the time series at $d$ timepoints, where $d$ is smaller than the given time sample. We assume timepoints $t = \frac{1}{d}, \frac{2}{d}, \ldots, \frac{d}{d}$. Let $\boldsymbol{t} = \left\{\frac{1}{d}, \frac{2}{d}, \ldots, 1\right\}$. So

$$X_i = \left(f_i\left(\frac{1}{d}\right), f_i\left(\frac{2}{d}\right), \ldots, f_i(1)\right)^T \in \mathcal{S}^d$$

We abstract from registration and filtering, and assume that appropriate subseries selection, detrending, and warping, rescaling and shifting into $\mathcal{S}^d \times [0,1]$, have taken place. Note that we impose relatively mild assumptions on the underlying functional form of our multivariate data. These assumptions are representative of the majority of Economic measurements, including functions that have jumps, like for instance the number of employees or the amount of capital, which could increase suddenly.

As we are dealing with multivariate time series, we are ultimately going to encounter curse of dimensionality issues. For this reason, and because of the discrete way digital computers work, we assume that we are evaluating the time series at $d$ more or less representative points. This way, we can choose the smallest possible $d$, such that the functions are sufficiently separated, when evaluated on $\boldsymbol{t}$ (in particular, we want to avoid the scenario in Figure 3). Section 5 contains some discussion about how to select a value for $d$.

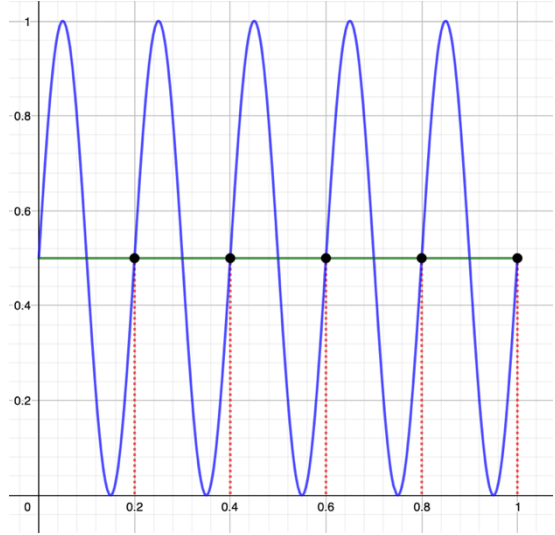

*Figure 3: For d=5, the blue and green function would be put in the same cluster*

We can always increase $d$ to closer approximate the true distance between functions. This is made precise by the following

**Lemma 4.0.1.** *For any $i \neq j$ we have $M_d \to M$ as $d \to \infty$, where $M_d = \sup_{t \in \mathbf{t}} \|f_i(t) - f_j(t)\|$, $M = \sup_{t \in [0,1]} \|f_i(t) - f_j(t)\|$. In particular, $\exists d' > d$, such that $d'' > d' \Rightarrow M_d \leq M_{d''} \leq M$.*

This brings us to estimation. Data may be incomplete, different time series may be sampled at different times, and we may not have datapoints at the times, we are evaluating the function for clustering ($t = \frac{1}{d}, \frac{2}{d}, \ldots, \frac{d}{d}$). In addition, most time series data are noisy, for instance due to measurement error (e.g., somebody mistyped a number in Excel). For these reasons, we need to estimate the functional form of the data before clustering. For series $X_i$, we assume we have a noisy time sample $Y_i$ of size $m_i > d$ in the interval $[0,1]$, at times $\frac{1}{m_i}, \frac{2}{m_i}, \ldots, 1$

Instead of $\frac{1}{m_i}, \frac{2}{m_i}, \ldots, 1$, for the result below it suffices to have sampling times, such that the samples grow in a sparse, uniform fashion across $[0,1]$, so that it will eventually get arbitrarily close to the points at which we evaluate the time series for clustering: the sampling frequency would need to increase in every subinterval.

So, with noise, we observe $y_i(t) = f_i(t) + \varepsilon_{it} \in \mathbb{R}^s$, maintaining the standard assumptions that

1. $\mathbb{E}\varepsilon_{it} = 0, \forall t, i$

2. $-\mathcal{E} \leq \varepsilon_{ikt} \leq \mathcal{E}$, for some $0 < \mathcal{E} < \infty$ where $\varepsilon_{ikt}$ are the component shocks
3. $\varepsilon_{it}$ are i.i.d., with independence across times $t$, components $k$ and time series $i$

Instead of 2., we could assume fast decreasing tails, i.e., $\mathbb{E}\varepsilon_{ikt}^4 < \infty$. To estimate $f_i$, we use a Nadaraya-Watson estimator

$$\hat{f}_i^{m_i}(t) = \sum_{j=1}^{m_i} w_j^{\gamma_{m_i}}(t) y_i(t_j)$$

$$w_j^{\gamma}(t) = \frac{K^{LE}\left(\frac{t-t_j}{\gamma}\right)}{\sum_{k=1}^{m_i} K^{LE}\left(\frac{t-t_k}{\gamma}\right)}$$

Where $\frac{0}{0} := 0$, $K^{LE}(x) = \mathbf{1}_{-1 \leq x \leq 0}(2 - 2x^2)$ is a left-side Epanechnikov kernel and $\gamma_m$ a choice of bandwidth, with $\gamma_m \to 0$, and $m\gamma_m \to \infty$, possible choices include $\gamma_m = m^{-\alpha}$, and $\gamma_m = \log m$. See Walk (2004) for a proof of strong universal consistency of the Nadaraya-Watson estimator. Our final estimated time series samples for clustering are then

$$\hat{X}_i = \left(\hat{f}_i^{m_i}\left(\frac{1}{d}\right), \hat{f}_i^{m_i}\left(\frac{2}{d}\right), \ldots \hat{f}_i^{m_i}(1)\right)^T$$

How to cluster this? For this, note that while comparing time series in $\mathcal{S}^d$ may not be clearly defined, we can flatten the series and every multivariate time series is then uniquely identified by a point in $[0,1]^{sd}$. We let

$$x_i = \left(f_i\left(\frac{1}{d}\right)^T, f_i\left(\frac{2}{d}\right)^T, \ldots, f_i(1)^T\right)^T \text{ and } \hat{x}_i = \left(\hat{f}_i^{m_i}\left(\frac{1}{d}\right)^T, \hat{f}_i^{m_i}\left(\frac{2}{d}\right)^T, \ldots, \hat{f}_i^{m_i}(1)^T\right)^T$$

the flattened versions of $X_i$ and $\hat{X}_i$, respectively.

We then construct the probability space $([0,1]^{sd}, \mathcal{F}, P)$, $\mathcal{F} = \mathcal{B}([0,1]^{sd})$ on all flattened time series $x_i$ as follows. We assume that each $X_i$ belongs to one of $C$ detectable clusters, where $C$ is unknown. Here, $C$ depends on $d$, and it is not necessarily the number of classes different functions $f_i$ (see Figure 3), but as $d \to \infty$, converges to the true number of clusters. Each detectable cluster is characterised by an unknown population frequency $\pi_1, \pi_2, \ldots, \pi_C$, where $\sum_{i=1}^{C} \pi_i = 1$ and a Lebesgue probability density $p_i$ fully supported on a compact subset of $[0,1]^{sd}$. The densities $p_i$ each define a measure $P_i(A) = \int_A p_i(x)\mu(dx)$, $A \in \mathcal{F}$ and the population mixture model is then $P = \sum_{i=1}^{C} \pi_i P_i$ characterised by the Lebesgue density $p(x) = \sum_{i=1}^{C} \pi_i p_i(x)$. We assume that each $p_i$ is $l_p$-Lipschitz, so

$$|p_i(x) - p_i(y)| \leq l_p|x - y|,$$

we assume unimodality, i.e.,

$$\text{for } a \in [0, \|p_i\|_\infty], \{x \in \mathbb{R}^{sd} : p_i \geq a\} \text{ is a convex set},$$

and further $x \notin [0,1]^{sd} \Rightarrow p_i(x) = 0$. We formalise the notion of detectability as follows. We assume that for each $i$ there exists a compact subset $S_i \subset [0,1]^{sd}$, such that

1. $S_i \cap S_j = \emptyset$, for $i \neq j$,
2. $x \in S_i \Rightarrow p_i(x) \geq 0$, and
3. $\forall x \notin S_i \Rightarrow p_i(x) < \lambda_*$, where $\lambda_*$ is unknown.

This way, the level set $\{p \geq \lambda_*\}$ has exactly $C$ connected components. Our assumptions rule out cusps and bridges and will allow us to derive well-behaved concentration bounds. They are implicit, i.e., for actual clustering purposes, we never have to worry about the $p_i$s, $S_i$s or $\lambda_*$ and it suffices to estimate the level sets of $p$.

By way of justification of these assumptions, we note the following. Every time series is theoretically defined at every point in the interval $[0,1]$. However, for clustering, a computer will necessarily only compare the data at a finite number of points, in our case, on $t$. For a given $d$, the underlying "distribution" of functions will induce a distribution on $[0,1]^{sd}$ and as we add more time series, we make i.i.d. draws from $p$. We can only cluster the data we have, and therefore it makes sense to check that our clustering algorithm is consistent, given values on $t$. We assume that functions of a similar "type" are close everywhere on $[0,1]$. This is formalised by the assumption that every time series is most likely to belong to only one subset $S_i$ of $[0,1]^{sd}$.

These assumptions match those made in the majority of the literature on clustering. Indeed, the clustering literature often additionally assumes the stronger condition of Hölder-continuous densities (see Wang et al., 2019; Sriperumbudur & Steinwart, 2012). We note that Lipschitz-continuity is not a necessary assumption, and we could make assumptions on the noisiness of the kernel density estimates as is done in Rinaldo and Wasserman (2011: Assumption C1), to allow for cases, where the densities $p_i$ are not defined.

We proceed in two steps. First, we establish the consistency and rates of DBSCAN for clustering the time series, assuming we observe the values the function takes exactly. This is in essence an application of kernel density estimation studied in Rinaldo and Wasserman (2010). In the second step, we introduce noise in the observed sample to the estimation problem and using the Nadaraya-Watson estimator, establish consistency.

## 5. Results

### 5.1 DBSCAN and density-based clustering

We note that some of the results in this section follow as an application of established results from the literature to our case. We use a simplified version of the DBSCAN algorithm for our case here, which is also used in Sriperumbudur and Steinwart (2012) and Wang *et al.* (2019). It goes as follows.

**Definition 5.1.1.** *(DBSCAN algorithm) For given k and δ.*

1. *For the sample of estimated time series $\widehat{D} = \{\hat{x}_1, \hat{x}_2, \ldots, \hat{x}_n\}$, identify $\mathcal{L}(k, \delta) = \{x \in \widehat{D} : |B(x, \delta) \cap \widehat{D}| \geq k\}$, i.e., sample points for which the closed balls $B(\hat{x}_i, \delta)$ contain at least k other samples.*
2. *The set of estimated clusters $\mathcal{T}(k, \delta)$ is then the set of δ-connected components of $\mathcal{L}(k, \delta)$.*

DBSCAN can be implemented to achieve this in sub $O(n^2)$ run time. For different values of $k$, $\mathcal{T}(k, \delta)$ is then an estimator of the cluster tree. To see this, note that for $k_1 > k_2$

$$\mathcal{L}(k_1, \delta) \subseteq \mathcal{L}(k_2, \delta)$$

So, $\mathcal{T}(k, \delta)$ necessarily exhibits the hierarchical characteristics of a cluster tree. In particular, as DBSCAN essentially uses a kernel density estimator $\hat{p}_\delta(x)$ with spherical kernel $K(x) = \mathbf{1}_{x \in B_{sd}(0,1)}$ (Sriperumbudur & Steinwart, 2012), we can relate $\mathcal{L}(k, \delta)$ to an estimate of the level sets of $\hat{p}_\delta$. $\delta \equiv \delta_n \geq 0$ is a bandwidth with $n\delta_n^{sd} \to \infty$.

For now, assume that we are given the actual time series $D = \{x_1, x_2, \ldots, x_n\}$. We will reintroduce estimation in a later section. We have

$$\hat{p}_\delta(x) = \frac{1}{n\delta^{sd} v_{sd}} \sum_{i=1}^{n} K\left(\frac{x - x_i}{\delta}\right)$$

Note that $\hat{p}_\delta(x)$ is a valid Lebesgue density, i.e.,

$$\int_{[-\delta, 1+\delta]^{sd}} \hat{p}_\delta(x)\mu(dx) = 1, \hat{p}_\delta(x) \geq 0$$

Furthermore, we have

$$p_\delta(x) \equiv \mathbb{E}[\hat{p}_\delta(x)] = \frac{P(B(x,\delta))}{\delta^{sd} v_{sd}}$$

And it can be shown that $\|\hat{p}_\delta - p_\delta\|_\infty \to 0, a.s.$ (see further below). $p_\delta(x)$ is what is referred to as the *mollified density* in Rinaldo and Wasserman (2011). In theory, mollification allows us to deal with cases, where the density is not defined. The mollified density defines a probability measure $P_\delta$ on $[-\delta, 1+\delta]^{sd}$

$$P_\delta(A) = \int_A p_\delta(x)\mu(dx)$$

$P_\delta$ and $p_\delta$ are approximations of $P$ and $p$ in the sense below, and in fact can be employed in a case where $p$ is not defined.

**Lemma 5.1.1.** *(Rinaldo & Wasserman, 2011) P-almost everywhere, $P_\delta$ converges weakly to P and $\lim_{\delta \to 0} p_\delta(x) = p(x)$*

Then, consider the estimator of the level set

$$\hat{L}(\lambda) = \{\hat{p}_\delta \geq \lambda\} = \cup_{x \in \{\hat{p}_\delta \geq \lambda\} \cap D} B(x, \delta)$$

This estimator is also studied in other places in the literature on kernel density estimation (see Devroye & Wise, 1980; Cuevas & Rodríguez-Casal, 2004; Wang *et al.*, 2019) and is in fact a consistent estimator of the level set $L(\lambda)$. Let $\lambda(k) = \frac{k}{n\delta^{sd} v_{sd}}$. From the definition $\hat{p}_\delta$ we directly have $x \in \mathcal{L}(k, \delta) \Rightarrow x \in \hat{L}(\lambda(k))$, since $|B(x,\delta) \cap D| \geq k \Leftrightarrow \sum_{i=1}^n K\left(\frac{x-x_i}{\delta}\right) \geq k$. Further, let $\hat{\mathcal{T}}(\lambda)$ be the set of connected components of $\hat{L}(\lambda)$, then

**Observation 5.1.1.**

$$A \in \hat{\mathcal{T}}(\lambda(k)) \Rightarrow A \cap D \in \mathcal{T}(k, \delta)$$

*and*

$A \in \mathcal{T}(k, \delta) \Rightarrow \exists \text{unique } A' \in \hat{\mathcal{T}}(\lambda(k))$, *such that* $A \subset A'$. *In particular,* $A' = \cup_{x \in A} B(x, \delta)$.

This means that $\hat{\mathcal{T}}(\lambda(k))$ and $\mathcal{T}(k, \delta)$ contain the same information and we can use the results on the consistency of $\hat{L}(\lambda)$ to prove the consistency of $\mathcal{T}(k, \delta)$. The trick here and indeed the beauty of DBSCAN, is that by computing $\mathcal{T}(k, \delta)$ instead of $\hat{\mathcal{T}}(\lambda)$, *we never have to estimate the density*.

## 5.2 Hartigan consistency of DBSCAN

To show that $\hat{L}(\lambda)$ is consistent, we need to show that the event $\|\hat{p}_\delta - p\|_\infty < \varepsilon$ holds with high probability as we increase $n$. Applying the triangle inequality, we have

$$\|\hat{p}_\delta - p\|_\infty \leq \|\hat{p}_\delta - p_\delta\|_\infty + \|p_\delta - p\|_\infty$$

By Devroye and Lugosi (2001: Corollary 4.2), we know that the collection of Euclidian balls on a compact set (in our case $[0,1]^{sd}$), has a finite Vapnik–Červonenkis dimension. Then, the concentration bound on the KDE is given by,

**Lemma 5.2.1.** $\mathbb{P}\left(\|\hat{p}_\delta - p_\delta\|_\infty < c\sqrt{\frac{\tau + \log 1/\delta}{n\delta^{sd}}}\right) \geq 1 - 2e^{-\tau}$, *for* $n\delta^{sd} \geq 1$, $\tau > 0$, *for some* $c_1 > 0$, *depending on* $\|p\|_\infty$, *s and d, independent of n and* $\tau$.

This result is originally due to Giné and Guillou (2002), and can be found in this form in Kim *et al.* (2019: Corollary 13). Bounds on $\|\hat{p}_\delta - p_\delta\|_\infty$ of the same form can also be found in Sriperumbudur and Steinwart (2012: Theorem 3.1) and Jiang (2017). Further, by the Lipschitz property,

**Lemma 5.2.2.** $\|p_\delta - p\|_\infty < c_2 \delta$ *where* $c_2$ *depends only on the Lipschitz constant of* $p$.

Combining Lemmata **5.2.1** and **5.2.2**, we have that with probability at least $1 - 2e^{-\tau}$,

$$\|\hat{p}_\delta - p\|_\infty \leq c_1 \sqrt{\frac{\tau + \log 1/\delta}{n\delta^{sd}}} + c_2 \delta$$

Then, letting $\tau = \log n$ and $\delta_n \coloneqq \left[z_1 \left(\frac{\log n}{n}\right)^{\frac{1}{2+sd}}, z_2 \left(\frac{\log n}{n}\right)^{\frac{1}{2+sd}}\right]$ for $z_2 \geq z_1 > 0$, we get $c_1 \sqrt{\frac{\tau + \log 1/\delta}{n\delta^{sd}}} + c_2 \delta \leq c' \left(\frac{\log n}{n}\right)^{\frac{1}{2+sd}}$, for an appropriate constant $c'$, so we have a final rate of $\left(\frac{\log n}{n}\right)^{\frac{1}{2+sd}}$, and

**Theorem 5.2.1.** *Let* $\delta_n \coloneqq \left(\frac{\log n}{n}\right)^{\frac{1}{2+sd}}$, *then for an appropriate constant* $c'$ *independent of n, with probability at least* $1 - \frac{2}{n}$

$$\|\hat{p}_{\delta_n} - p\|_\infty \leq c' \left(\frac{\log n}{n}\right)^{\frac{1}{2+sd}}$$

The following corollary is then of interest.

**Corollary 5.2.1.** *Let* $\delta_n \coloneqq \left(\frac{\log n}{n}\right)^{\frac{1}{2+sd}}$ *then,*

$$\mathbb{P}\left\{L\left(\lambda + c'\left(\frac{\log n}{n}\right)^{\frac{1}{2+sd}}\right) \subset \hat{L}(\lambda) \subset L\left(\lambda - c'\left(\frac{\log n}{n}\right)^{\frac{1}{2+sd}}\right)\right\} \geq 1 - \frac{2}{n}$$

Which shows consistency of $\hat{L}(\lambda)$.

Finally, we can show Hartigan consistency.

**Theorem 5.2.2.** $\hat{T}(\lambda)$ *is* Hartigan consistent. *In particular, given a sample $D_n$ of size $n$ for any disjoint subsets $A, A'$ of $[0,1]^{sd}$, such that $A, A' \in L(\lambda)$ are in disconnected components of $L(\lambda)$ for some $\lambda$, we have $P(A_n \cap A'_n = \emptyset) \geq 1 - \frac{2}{n}$, where $A_n, A'_n$ are the smallest clusters in $\hat{T}$ that contain $A \cap D_n$ and $A' \cap D_n$, respectively.*

Then

**Corollary 5.2.2.** $\mathcal{T}(k(\lambda), \delta_n)$ *with* $k(\lambda) = \lambda n \delta_n^{sd} v_{sd}$ *and* $\delta_n \asymp \left(\frac{\log n}{n}\right)^{\frac{1}{2+sd}}$ *is a Hartigan consistent estimator of the cluster tree in the sense of Observation* **5.2.1.**

Thus, we have established that DBSCAN is consistent for the flattened problem.

## 5.3 Introducing noise and estimation

In this section, we will include the estimation step. In our setup, and indeed in most practical clustering scenarios, we first estimate the $x_i$s using kernel smoothing to get estimates $\hat{x}_i$. As such, our KDE is $\hat{\tilde{p}}_{\delta,\gamma}(x) = \frac{1}{n\delta^{sd}v_{sd}} \sum_{i=1}^{n} K\left(\frac{x-\hat{x}_i}{\delta}\right)$. To achieve consistency, then, we need the event

$$E_{\delta,\gamma,\varepsilon} = \left\{\left\|\hat{\tilde{p}}_{\delta,\gamma} - p\right\|_{\infty} < \varepsilon\right\}$$

to hold with high probability as $m, m\gamma_m \to \infty$, $n, n\delta_n^{sd} \to \infty$, $\gamma_m, \delta_n \to 0$, where $m \equiv \min_i m_i$. The aim of this section is to derive consistency of $\hat{\tilde{p}}_{\delta,\gamma}$. For this, we derive upper bounds on the convergence rates, that are not necessarily the sharpest bounds possible, but which suffice to demonstrate convergence.

An application of the triangle inequality yields,

$$\left\|\hat{\tilde{p}}_{\delta,\gamma} - p\right\|_{\infty} \leq \left\|\hat{\tilde{p}}_{\delta,\gamma} - \mathbb{E}\hat{\tilde{p}}_{\delta,\gamma}\right\|_{\infty} + \left\|\mathbb{E}\hat{\tilde{p}}_{\delta,\gamma} - p_{\delta}\right\|_{\infty} + \left\|p_{\delta} - p\right\|_{\infty}$$

We have already established bounds for $\|p_\delta - p\|_\infty$ in the previous section and $\left\|\hat{p}_{\delta,\gamma} - \mathbb{E}\hat{p}_{\delta,\gamma}\right\|_\infty \to 0$ can be guaranteed by a similar argument as before. To bound the bias $\left\|\mathbb{E}\hat{p}_{\delta,\gamma} - p_\delta\right\|_\infty$, we first establish the behaviour of $\hat{p}_{\delta,\gamma} - \hat{p}_\delta$.

We begin by establishing the convergence of $\|x_i - \hat{x}_i\|_\infty$. The functions $f_i$ are estimated using a Nadaraya-Watson kernel regression. Specifically, we use a left-side kernel. This is because we only assume right-continuity and the $f_i$s may contain jumps. Estimating $f_i$ using a symmetric kernel, the Nadaraya-Watson estimator would fail to consistently estimate $f_i$ at jump points, as the kernel will always contain samples to the right of the jump for every bandwidth size, meaning that the bias persists as $m \to \infty$, and that we could not establish uniform bounds. This problem is illustrated in Figure 4. A left-side kernel does not have this problem. We note that for low values of $m_i$, the left-side kernel will still have additional bias at points after jumps. The difference is that when the sample size is large enough, this "jump-bias" completely vanishes.

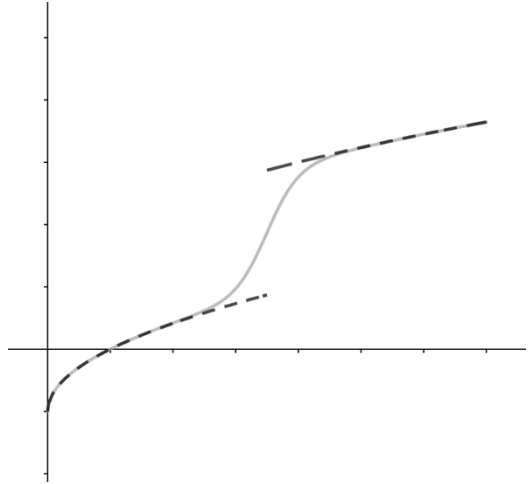

*Figure 4: Right continuous function (dashed line) and Nadaraya-Watson estimator (in gray) using a symmetric kernel. The estimate at the jump will always lie half-way between the right and left limit at the jump.*

Given this, we first establish

**Theorem 5.3.1.** $\forall d \in \mathbb{N}$ as $m \to \infty$, $\sup\limits_{t \in t}\left\|\hat{f}_i^{m_i}(t) - f_i(t)\right\|_\infty \to 0, a.s.$ In particular, for $m$ large enough and with probability at least $1 - 2e^{-\tau}$

$$\sup_{t \in t}\left\|\hat{f}_i^{m_i}(t) - f_i(t)\right\|_\infty \leq c'\sqrt{\frac{\tau - \log sd}{m\gamma_m}} + l_f \gamma_m$$

*Where $c'$ is an appropriate constant and $l_f$ is the Lipschitz constant for $f_i$.*

Then, the following corollary is immediate.

**Corollary 5.3.1.** $\|x_i - \hat{x}_i\|_\infty \leq c'\sqrt{\frac{\tau - \log sd}{m\gamma_m}} + l_f\gamma_m$ and $\sup_i \|x_i - \hat{x}_i\|_\infty \leq c'\sqrt{\frac{\tau - \log sdn}{m\gamma_m}} + l_f\gamma_m$ with probability at least $1 - 2e^{-\tau}$ in each case.

Then let $\tau = \log m$ and $\gamma_m \asymp \left(\frac{\log m}{m}\right)^{\frac{1}{3}}$. This means that for $n$ given, $\sup_i \|x_i - \hat{x}_i\|_\infty$ decreases like $\left(\frac{\log m}{m}\right)^{\frac{1}{3}}$. It is important to note that the probability of the event $\|x_i - \hat{x}_i\|_\infty < \varepsilon$ depends only on $m$ and is independent of $x_i$.

A crude consistency argument then goes as follows. As $\sup_i \|x_i - \hat{x}_i\|_\infty \to 0$, we have $P_\delta(\{x \in [-\delta, 1+\delta]^{sd}: |\hat{\hat{p}}_{\delta,\gamma}(x) - \hat{p}_\delta(x)| = 0\}) \to 1$, where $P_\delta(A) = \int_{x \in A} p_\delta(x)\mu(dx)$. That is, the set on which the two KDEs coincide, converges to the support of $\hat{p}_\delta$. In the limit, $\hat{p}_\delta(x) = \hat{\hat{p}}_{\delta,\gamma}(x)$, $P_\delta$-almost everywhere. To show this, we note,

$$|\hat{\hat{p}}_{\delta,\gamma}(x) - \hat{p}_\delta(x)| = \left|\frac{1}{n\delta^{sd}v_{sd}}\sum_{i=1}^{n}\left(\mathbf{1}_{x\in B(\hat{x}_i,\delta)} - \mathbf{1}_{x\in B(x_i,\delta)}\right)\right|$$

$$\leq \frac{1}{n\delta^{sd}v_{sd}}\sum_{i=1}^{n}\left|\mathbf{1}_{x\in B(\hat{x}_i,\delta)} - \mathbf{1}_{x\in B(x_i,\delta)}\right|$$

Where specifically, we have

$$\left|\mathbf{1}_{x\in B(\hat{x}_i,\delta)} - \mathbf{1}_{x\in B(x_i,\delta)}\right| = \left|\mathbf{1}_{x\in B(\hat{x}_i,\delta)\setminus B(x_i,\delta)} - \mathbf{1}_{x\in B(x_i,\delta)\setminus B(\hat{x}_i,\delta)}\right| = \mathbf{1}_{x\in B(x_i,\delta)\triangle B(\hat{x}_i,\delta)}$$

To simplify notation, let $B_i \equiv B(x_i, \delta)$ and $\hat{B}_i \equiv B(\hat{x}_i, \delta)$. Then

**Lemma 5.3.1.** *For a time sample of size $m$ and given $\delta$, let the $E_{m,\varepsilon} = \left\{\sup_i \|x_i - \hat{x}_i\|_\infty \leq \frac{\varepsilon}{\sqrt{sd}}\right\}$ hold, then $\mu(B_i \triangle \hat{B}_i) \leq 2^{sd}\delta^{sd-1}\varepsilon$*

**Remark 5.3.1.** *The inequality in Lemma 5.3.1 can be sharpened to*

$$\mu(B_i \triangle \hat{B}_i) = 2\delta^{sd}v_{sd}\left(1 - I_{1-\frac{\varepsilon^2}{4\delta^2}}\left(\frac{sd+1}{2}, \frac{1}{2}\right)\right)$$

*For some $c$ and where $I_x(a, b)$ is the regularised incomplete beta function.*

In particular, this means that the volume of $\cup_i B_i \triangle \hat{B}_i$ decreases at the same rate as $\|x_i - \hat{x}_i\|_\infty$ and hence $P_\delta(\{x \in [-\delta, 1+\delta]^{sd}: |\hat{\hat{p}}_{\delta,\gamma}(x) - \hat{p}_\delta(x)| > 0\})$ will decrease at a rate proportional to that of $\|x_i - \hat{x}_i\|_\infty$, because we have

$$\{x \in [-\delta, 1+\delta]^{sd}: |\hat{\hat{p}}_{\delta,\gamma}(x) - \hat{p}_\delta(x)| = 0\} = \{x \in [-\delta, 1+\delta]^{sd}: x \in (\cup_i B_i \triangle \hat{B}_i)^c\}$$

By subadditivity, we have $\mu(\cup_i B_i \triangle \hat{B}_i) \leq 2^{sd} n \delta_n^{sd-1} \varepsilon$, so if we simultaneously increase $n$, we need to do so sufficiently slowly, such that $\frac{\varepsilon_m}{n \delta_n^{sd-1}} \to 0$ as $m, n \to \infty$.

Note that the "random variable" is $B_i \triangle \hat{B}_i$, which depends on the joint distribution of $x_i$ and $\hat{x}_i$, conditional on the event $E_{m,\varepsilon}$. We can show,

**Lemma 5.3.2.** *For given $\gamma_m$, $\delta$, $x$ and $\varepsilon \leq \delta$, on the event $E_{m,\varepsilon} = \{\sup_i \|x_i - \hat{x}_i\|_\infty \leq \frac{\varepsilon}{\sqrt{sd}}\}$, we have*

$$\mathbb{E}|\hat{\hat{p}}_{\delta,\gamma}(x) - \hat{p}_\delta(x)| \leq \frac{z'\varepsilon}{\delta}$$

*For $z' > 0$, where $z'$ depends on $\|p\|_\infty$, and $v_{sd}$. Furthermore, with probability at least $1 - \tau$*

$$|\hat{\hat{p}}_{\delta,\gamma}(x) - \hat{p}_\delta(x)| \leq \frac{z'\varepsilon}{\delta\tau}$$

This establishes pointwise convergence $\hat{\hat{p}}_{\delta,\gamma}(x) \to^p \hat{p}_\delta(x)$, take $\gamma_m \asymp \left(\frac{\log m}{m}\right)^{\frac{1}{3}}$, $\varepsilon_m = c'' \left(\frac{\log m}{m}\right)^{\frac{1}{3}}$ and $\tau = \sqrt{\left(\frac{\log m}{m}\right)^{\frac{1}{3}}}$, then for a given $n$, $|\hat{\hat{p}}_{\delta,\gamma}(x) - \hat{p}_\delta(x)|$ decreases like $\sqrt{\left(\frac{\log m}{m}\right)^{\frac{1}{3}}}$.

As a direct corollary, we can establish the desired bound for $\|\mathbb{E}\hat{\hat{p}}_{\delta,\gamma} - p_\delta\|_\infty$.

**Theorem 5.3.2.** *With probability at least $1 - 2e^{-\tau}$ and for some $z'' > 0$,*

$$\|\mathbb{E}\hat{\hat{p}}_{\delta,\gamma} - p_\delta\|_\infty \leq z'' \sqrt{\frac{\tau - \log n}{m \gamma_m \delta^2}}$$

Which gives us consistency of $\hat{\hat{p}}_{\delta,\gamma}$.

# 6. Practical considerations for the choice of $d$

A very important choice to make in the setup of this paper is the choice of $d$. As we increase $d$, we will increase to accuracy with which we classify different functions. However, we simultaneously increase the volume of the confidence set of the level set estimator unless we increase $n$ sufficiently.

The rate of almost sure convergence is given by $r_n^d = \left(\frac{\log n}{n}\right)^{\frac{1}{2+sd}}$, so for $n > 1$

$\frac{\partial r_n^d}{\partial d} = \frac{-s}{(2+sd)^2} \log\left(\frac{\log n}{n}\right) \left(\frac{\log n}{n}\right)^{\frac{1}{2+sd}} > 0$, which means, as expected, the rate becomes less sharp as we increase $d$, which is the curse of dimensionality.

Let $\Delta$ be the compensating differential, such that, $r_{n+\Delta}^{d+1} = r_n^d$. Solving for $\Delta$ yields

$$\Delta \approx 0.74 \left(\frac{n}{\log n}\right)^{2+\frac{2s}{2+sd}} - n$$

This grows at a super-linear rate. As we increase $d$, to regain the same performance, we essentially have to square the sample size, save for log factors. This highlights, why the choice of $d$ should be conservative. It is only worth increasing $d$, if it substantially changes the composition of the clusters. As there is no correct $d$, however, it is best to report results for a range of $d$, as is common in the clustering literature, e.g., with the choice of $k$ in $k$-means.

# 7. Conclusion and Discussion

This paper has established consistency results and rates for clustering flattened multivariate data using DBSCAN without any dimensionality reduction. As DBSCAN is a computationally efficient clustering algorithm, this suggests that using it as a subroutine in the way studied in this paper may be highly useful for empirical research. As is usual with non-parametric methods, however, we are dealing with exceedingly slow convergence rates, slower than $\frac{1}{\sqrt{n}}$ and with additional log factors, so this way of clustering data is best suited for large samples of high frequency data. This is not a problem for most Economic data, as high-frequency data is either readily available (e.g., financial data) or may be reconstructed (e.g., the number of employees is the same at times $t$ and $t + \varepsilon$).

There are many areas of inquiry that could be investigated in future research. Notably,

1. Further research should be concerned with establishing sharper bounds, and the choice of suitable bandwidths given $n$ and $m$.
2. An area of particular interest would be to use the results of this paper as a reference point for quantifying the bias that results by first employing dimensionality reduction, e.g., in the form of FPCA.
3. This paper has fully abstracted from filtering and registration. Further research could investigate the consistency implications of various filtering algorithms.

As a final epistemological remark, note that the notion of consistency, in this paper and the wider literature, while mathematically elegant, may be a meaningless concept in terms of the real problem we are solving. For most *things*, there is no such thing as the infinite population, and perhaps we had better investigate the mechanisms, by which a finite sample is chosen out of a finite population.

**Proofs**

**Lemma 4.0.1.** Trivial: [0,1] is compact, so $\exists t^*$, such that $M = \|f_i(t^*) - f_j(t^*)\|$. $f_i, f_j$ are right continuous so $\Delta(t) = \|f_i(t) - f_j(t)\|$ is also right continuous. Let $t_<(d) = \sup\{t \in t : t \leq t^*\}$. As $t^* - t_<(d) < \frac{1}{d}$, we have that $t_<(d)$ is increasing. By right-continuity of $\Delta(t)$, $\Delta(t_<(d)) \to \Delta(t^*) = M$, in fact since $\Delta(t_<(d)) \leq M$, $\Delta(t_<(d))$ is also increasing. Since $M \geq M_d \geq \Delta(t_<(d))$, we have $M_d \to M$ by the Sandwich theorem. The last claim follows from the fact that $\Delta(t_<(d))$ is increasing.

**Lemma 5.1.1.** A more general version of this proof, including for cases where the density is not defined, can be found in Rinaldo and Wasserman (2011). The weak convergence is immediate from the fact that $P$ is Radon. The second property can be demonstrated as follows. We have

$$p_\delta(x) = \frac{P(B(x,\delta))}{\delta^{sd} v_{sd}}$$

Note that

$$\delta^{sd} v_{sd} \inf_{y \in B(x,\delta)} p(y) \leq P(B(x,\delta)) \leq \delta^{sd} v_{sd} \sup_{y \in B(x,\delta)} p(y)$$

so $\inf_{y \in B(x,\delta)} p(y) \leq p_\delta(x) \leq \sup_{y \in B(x,\delta)} p(y)$.

Then by continuity, $\inf_{y \in B(x,\delta)} p(y) \to p(x)$ and $\sup_{y \in B(x,\delta)} p(y) \to p(x)$, and it is immediate that $p_\delta(x) \to p(x)$, as required.

**Lemma 5.2.2.** Same reasoning as in Lemma **5.1.1.**, and using the Lipschitz-ness of $p$. ∎

**Corollary 5.2.1.** We begin by proving the following lemma

*Lemma.* On the event $\|\hat{p}_\delta - p\| \leq \varepsilon$, we have $L(\lambda + \varepsilon) \subset \hat{L}(\lambda) \subset L(\lambda - \varepsilon)$

*Proof.* Consider the level set $\hat{L}(\lambda) = \{\hat{p}_\delta \geq \lambda\} = \cup_{x \in \{\hat{p}_\delta \geq \lambda\} \cap D} B(x,\delta)$

On the event $\|\hat{p}_\delta - p\| \leq \varepsilon$, we have

$$|\hat{p}_\delta(x) - p(x)| \leq \|\hat{p}_\delta - p\|_\infty \Leftrightarrow p(x) - \varepsilon \leq \hat{p}_\delta(x) \leq p(x) + \varepsilon$$

Take $x$, such that $\hat{p}_\delta(x) \geq \lambda$, by $\hat{p}_\delta(x) \leq p(x) + \varepsilon$, we get

$$\hat{p}_\delta(x) \geq \lambda \Rightarrow p(x) \geq \lambda - \varepsilon$$

Then clearly the set $\{x: \hat{p}_\delta(x) \geq \lambda\} \subseteq \{x: p(x) \geq \lambda - \varepsilon\}$, so $\hat{L}(\lambda) \subseteq L(\lambda - \varepsilon)$

By the same argument, take $x$ such that $p(x) \geq \lambda + \varepsilon$, then by $p(x) - \varepsilon \leq \hat{p}_\delta(x)$, we get $\hat{p}_\delta(x) + \varepsilon \geq p(x) \geq \lambda + \varepsilon \Rightarrow \hat{p}_\delta(x) \geq \lambda$ and thus $\hat{L}(\lambda) \supseteq L(\lambda + \varepsilon)$.

Together, this gives $L(\lambda + \varepsilon) \subseteq \hat{L}(\lambda) \subseteq L(\lambda - \varepsilon)$ as required. □

Then, the desired concentration bounds for $\hat{L}(\lambda)$ follows directly from the concentration bound in Theorem **5.2.1.**

**Theorem 5.2.2.**

Let $A, A'$ be disjoint subsets of $[0,1]^{sd}$. We have $\lambda = \inf_{x \in A \cup A'} p(x)$ and

$\underline{\lambda} = \sup_{\{P\ path\ from\ A\ to\ A'\ on\ [0,1]^{sd}\}} \inf_{x \in P} p(x)$. Then $A, A' \subset L(\lambda) \subset L(\underline{\lambda})$. In particular, $L(\underline{\lambda})$ is the highest level set, in which $A$ and $A'$ can be connected by some path.

Let the event $\|\hat{p}_\delta - p\|_\infty \leq \varepsilon_n$ hold. Applying the auxiliary Lemma from Corollary **5.2.1**, when $\lambda - \underline{\lambda} > 2\varepsilon_n$, then $A, A' \subset L(\lambda) \subset \hat{L}(\lambda - \varepsilon_n)$. However,

$$\hat{L}(\lambda - \varepsilon_n) \subset L(\lambda - 2\varepsilon_n) \subset L(\underline{\lambda})$$

And since $\lambda - 2\varepsilon_n > \underline{\lambda}$, it is immediate from the definition of $\underline{\lambda}$, that there is no path in $L(\lambda - 2\varepsilon_n)$ and hence in $\hat{L}(\lambda - \varepsilon_n)$ connecting $A$ to $A'$. This means that $A$ and $A'$ lie in disjoint connected components of $\hat{L}(\lambda - \varepsilon_n)$.

Finally, let $\delta_n \coloneqq \left(\frac{\log n}{n}\right)^{\frac{1}{2+sd}}$ and $\varepsilon_n = c'\left(\frac{\log n}{n}\right)^{\frac{1}{2+sd}}$ for some appropriate $c' > 0$. Then, by Theorem **5.2.1**, we get

$$\mathbb{P}(A_n \cap A'_n = \emptyset) \geq 1 + \frac{2}{n} \to 1$$

As required. Further, as $n \to \infty$, $\varepsilon_n \to 0$ and this eventually holds for all disjoint sets $A, A'$. Thus, we have established Hartigan consistency.

**Theorem 5.3.1.** Note $\sup_{t \in t}\|\hat{f}_i^{m_i} - f_i\|_\infty \leq \sup_{t \in t}\|f_i - \mathbb{E}\hat{f}_i^{m_i}\|_\infty + \sup_{t \in t}\|\hat{f}_i^{m_i} - \mathbb{E}\hat{f}_i^{m_i}\|_\infty$.

First consider $\left|f_i^k(t) - \mathbb{E}\left[(\hat{f}_i^{m_i})^k(t)\right]\right| = \left|\sum_{j=1}^{m_i}\left(f_i(t_j) - f_i(t)\right)w_j^{\gamma_{m_i}}(t)\right|$

$$= \left|\sum_{t-\gamma_{m_i} \leq t_j \leq 0}\left(f_i(t_j) - f_i(t)\right)w_j^{\gamma_{m_i}}(t)\right|$$

$$\leq \sum_{t-\gamma_{m_i} \leq t_j \leq 0}|f_i(t_j) - f_i(t)|w_j^{\gamma_{m_i}}(t) \leq l_f\gamma_{m_i} \leq l_f\gamma_m$$

Applying the Lipschitz property. As a direct consequence $\sup_{t \in [0,1]}\|f_i - \mathbb{E}\hat{f}_i^{m_i}\|_\infty \leq l_f\gamma_m$

Secondly, consider $\left|\mathbb{E}\left[(\hat{f}_i^{m_i})^k(t)\right] - (\hat{f}_i^{m_i})^k(t)\right| = \left|\sum_{j \leq m_i}\varepsilon_{ikt_j}w_j^{\gamma_{m_i}}(t)\right|$

Consider the set $\{j: t - \gamma_{m_i} \leq t_j \leq 0\} \equiv S_{m_i}$. This set contains between $\lfloor m_i\gamma_{m_i}\rfloor$ and $\lfloor m_i\gamma_{m_i}\rfloor - 1$ samples.

We have $-\mathcal{E} \leq \varepsilon_{ikt_j} \leq \mathcal{E}$. Also,

$$w_j^{\gamma_{m_i}}(t) \leq \frac{1 - (t - t_j)^2}{\sum_{q \in S_{m_i}}\left(1 - (t - t_q)^2\right)} \leq \frac{1}{\sum_{q \in S_{m_i}}(1 - m_i^{-2})} \leq \frac{1}{(\lfloor m_i\gamma_{m_i}\rfloor - 1)(1 - m_i^{-2})}$$

Therefore when $j \in S_{m_i}$,

$$-\frac{\mathcal{E}}{(\lfloor m_i\gamma_{m_i}\rfloor - 1)(1 - m_i^{-2})} \leq w_j^{\gamma_{m_i}}(t)\varepsilon_{ikt_j} \leq \frac{\mathcal{E}}{(\lfloor m_i\gamma_{m_i}\rfloor - 1)(1 - m_i^{-2})}$$

Then a direct application of Hoeffding's inequality yields

$$\mathbb{P}\left(\left|\mathbb{E}\left[(\hat{f}_i^{m_i})^k(t)\right] - (\hat{f}_i^{m_i})^k(t)\right| < \sqrt{2\tau\mathcal{E}^2 \sum_{j \in S_{m_i}} (\lfloor m_i\gamma_{m_i}\rfloor - 1)^{-2}(1 - m_i^{-2})^{-2}}\right)$$
$$\geq 1 - 2e^{-\tau}$$

Further simplification yields, assuming $m_i > 1$ and for an appropriate $c > 0$

$$\sum_{j \in S_{m_i}} (\lfloor m_i\gamma_{m_i}\rfloor - 1)^{-2}(1 - m_i^{-2})^{-2} \leq m_i\gamma_{m_i}(\lfloor m_i\gamma_{m_i}\rfloor - 1)^{-2}(1 - m_i^{-2})^{-2}$$

$$\leq 2m_i\gamma_{m_i}(\lfloor m_i\gamma_{m_i}\rfloor - 1)^{-2} \leq \frac{2c}{m_i\gamma_{m_i}} \leq \frac{2c}{m\gamma_m}$$

So, with $c' = 4c\mathcal{E}^2$

$$\mathbb{P}\left(\left|\mathbb{E}\left[(\hat{f}_i^{m_i})^k(t)\right] - (\hat{f}_i^{m_i})^k(t)\right| < c'\sqrt{\frac{\tau}{m\gamma_m}}\right) \geq 1 - 2e^{-\tau}$$

Using the union bound property/Bernoulli's inequality

$$\mathbb{P}\left(\sup_{t \in t}\|\hat{f}_i^{m_i}(t) - \mathbb{E}\hat{f}_i^{m_i}(t)\|_\infty < c'\sqrt{\frac{\tau - \log sd}{m\gamma_m}}\right) \geq 1 - 2e^{-\tau}$$

And we get the desired result.

Note that had we assumed instead that the component shocks have bounded fourth moment, an application of Kolmogorov's strong law of large numbers yields

$$P\left(\left|\mathbb{E}\left[(\hat{f}_i^{m_i})^k(t)\right] - (\hat{f}_i^{m_i})^k(t)\right| < \varepsilon\right) \geq 1 - O\left(\frac{1}{m_i^2\gamma_{m_i}^2}\right)$$

**Lemma 5.3.1.** The argument goes as follows, for a visual reference see Figure 5 below. $B(x_i, \delta)$ and $B(\hat{x}_i, \delta)$ are hyperspheres with overlap determined by $\|x_i - \hat{x}_i\|$. We have $\mu\big(B(x_i, \delta)\backslash B(\hat{x}_i, \delta)\big) = \mu\big(B(\hat{x}_i, \delta)\backslash B(x_i, \delta)\big)$ by definition. On any plane $P$ that contains $x_i$ and $\hat{x}_i$, the intersection of $B(x_i, \delta)$ and $B(\hat{x}_i, \delta)$ are two circles that intersect according to $\|x_i - \hat{x}_i\|$. These two circles intersect twice, say at points $A, B$, where $\|A - B\| \leq 2\delta$, the diameter of the circles. Then, if we trace the line segments $\overline{AB}$, this will form a hypersphere of dimension $sd - 1$. Then take the set $B(x_i, \delta)\backslash B(\hat{x}_i, \delta)$. This set forms an $sd$-dimensional "hyper-half-moon", which after a volume preserving transformation, can be completely contained in an $sd$-hypercylinder with radius $\delta$ and height $\|x_i - \hat{x}_i\|$. Then for the volume of the set $B(x_i, \delta) \triangle B(\hat{x}_i, \delta)$, we get

$$\mu\big(B(x_i, \delta) \triangle B(\hat{x}_i, \delta)\big) < 2^{sd} \delta^{sd-1} \|x_i - \hat{x}_i\|$$

On the event $E_{m,\varepsilon}$, we have $\|x_i - \hat{x}_i\| \leq \varepsilon \ \forall i.$, so we get the desired result,

$$\mu\big(B(x_i, \delta) \triangle B(\hat{x}_i, \delta)\big) < 2^{sd} \delta^{sd-1} \varepsilon$$

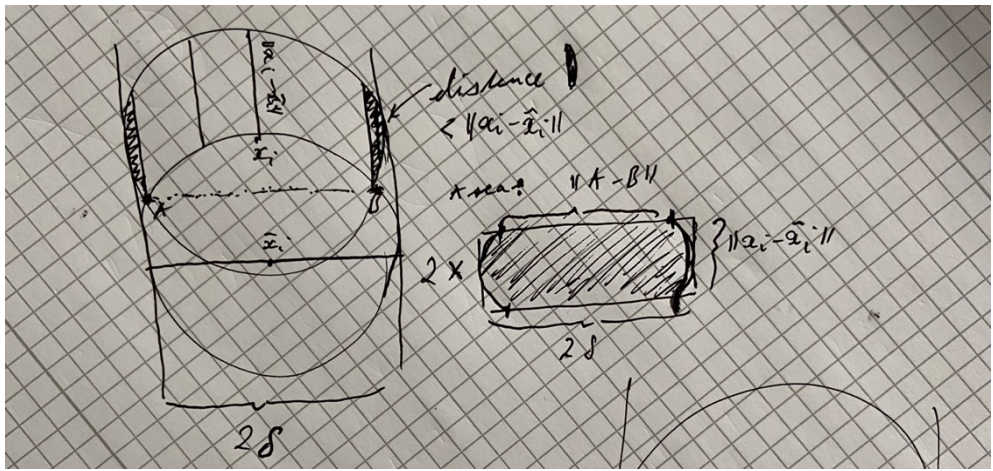

*Figure 5: Visual reference*

**Remark 5.3.1.**

We have $\mu\big(B(\hat{x}_i, \delta) \backslash B(x_i, \delta)\big) = \mu\big(B(x_i, \delta) \backslash B(\hat{x}_i, \delta)\big) = \delta^{sd} v_{sd} - \mu\big(B(x_i, \delta) \cap B(\hat{x}_i, \delta)\big)$, where

$$\mu\big(B(x_i, \delta) \cap B(\hat{x}_i, \delta)\big) = \delta^{sd} v_{sd} I_{1-\frac{\|x_i - \hat{x}_i\|^2}{4\delta^2}} \left(\frac{sd+1}{2}, \frac{1}{2}\right)$$

See Li (2011) for a concise derivation of the volume of the intersection of two $n$-balls. Then

$\mu\big(B(x_i, \delta) \backslash B(\hat{x}_i, \delta)\big) = \delta^{sd} v_{sd} \left(1 - I_{1-\frac{\|x_i - \hat{x}_i\|^2}{4\delta^2}} \left(\frac{sd+1}{2}, \frac{1}{2}\right)\right)$. The remark follows. This rate is slightly sharper than the one in Lemma **5.3.1.** because it additionally accounts for the decreasing radius of the $(sd - 1)$-hypersphere intersection between the balls.

**Lemma 5.3.2.**

First note $\mathbb{E}\big[|\mathbf{1}_{x \in B_i} - \mathbf{1}_{x \in \hat{B}_i}|\big] = \mathbb{E}[\mathbf{1}_{x \in B_i \triangle \hat{B}_i}]$

Let the event $\sup_i \|x_i - \hat{x}_i\|_\infty \leq \frac{\varepsilon}{\sqrt{sd}}$ hold and let $b_i \equiv \hat{x}_i - x_i$ be the bias, which is also a random variable. In particular, then $\hat{x}_i = x_i + b_i$. Consider the ball $B(x, \delta)$. We have $\mathbf{1}_{x \in B_i \triangle \hat{B}_i} = 1$, whenever the addition of the bias to $x_i$ changes, whether the point is in $B(x, \delta)$. So

$$\mathbf{1}_{x \in B_i \triangle \hat{B}_i} = \mathbf{1}_{\{x_i, b_i : x_i \in B(x,\delta) \text{ and } x_i + b_i \notin B(x,\delta)\} \cup \{x_i, b_i : x_i \notin B(x,\delta) \text{ and } x_i + b_i \in B(x,\delta)\}}$$

In order for the addition of $b_i$ to change the "status" of $x_i$, however, $x_i$ has to lie within $\|b_i\|$ of $\partial B(x, \delta)$, the boundary of the $\delta$-ball around $x$. On the event $\sup_i \|b_i\|_\infty \leq \frac{\varepsilon}{\sqrt{sd}}$, with $\varepsilon \leq \delta$, therefore

$$\{x_i, b_i : x_i \in B(x, \delta) \text{ and } x_i + b_i \notin B(x, \delta)\} \cup \{x_i, b_i : x_i \notin B(x, \delta) \text{ and } x_i + b_i \in B(x, \delta)\}$$
$$\subset \{x_i, b_i : \delta - \varepsilon \leq \|x_i - x\| \leq \delta + \varepsilon\}$$

Or

$$\{x_i, b_i : x \in B_i \triangle \hat{B}_i\} \subset \{x_i, b_i : \delta - \varepsilon \leq \|x_i - x\| \leq \delta + \varepsilon\}$$

The strict containment relation is due to the fact that while $x_i$ needs to lie within $\varepsilon$ of the boundary, not every $b_i$ will effect a shift across the boundary. Then, in particular,

$$\mathbb{E}[\mathbf{1}_{x \in B_i \triangle \hat{B}_i}] \leq \mathbb{P}(\{x_i, b_i : \delta - \varepsilon \leq \|x_i - x\| \leq \delta + \varepsilon\}) = \mathbb{P}(\{x_i : \delta - \varepsilon \leq \|x_i - x\| \leq \delta + \varepsilon\})$$
$$\leq \|p\|_\infty v_{sd}((\delta + \varepsilon)^{sd} - (\delta - \varepsilon)^{sd}) \leq z' \delta^{sd-1} v_{sd} \varepsilon$$

With $z' = 2^{sd} sd \|p\|_\infty$. Where the last inequality follows from the mean value theorem: $\exists x \in [\delta - \varepsilon, \delta + \varepsilon]$, such that

$$|(\delta + \varepsilon)^{sd} - (\delta - \varepsilon)^{sd}| = sd |x^{sd-1}| |\delta + \varepsilon - (\delta - \varepsilon)|$$

Let $\delta \geq \varepsilon$, then $|x^{sd-1}|$ can be bounded above by $(2\delta)^{sd-1}$, and we get the inequality.

Lastly, $\mathbb{E}|\hat{\hat{p}}_{\delta,\gamma}(x) - \hat{p}_\delta(x)| \leq \frac{1}{n \delta^{sd} v_{sd}} \sum_i \mathbb{E}[\mathbf{1}_{x \in B_i \triangle \hat{B}_i}] \leq \frac{z' \varepsilon}{\delta}$

Applying the Markov inequality, on the event $\sup_i \|x_i - \hat{x}_i\| \leq \varepsilon$ we get

$$\mathbb{P}\left(|\hat{\hat{p}}_{\delta,\gamma}(x) - \hat{p}_\delta(x)| < \frac{z' \varepsilon}{\tau \delta}\right) \geq 1 - \tau$$

And hence the desired result.

**Theorem 5.3.2.** On the event $E_{m,\varepsilon} = \left\{\sup_i \|x_i - \hat{x}_i\|_\infty \leq \frac{\varepsilon}{\sqrt{sd}}\right\}$, by Lemma **5.3.2.**, we always have

$$\mathbb{E}|\hat{\hat{p}}_{\delta,\gamma}(x) - \hat{p}_\delta(x)| \leq \frac{z'\varepsilon}{\delta} \; \forall x$$

Then $\left\|\mathbb{E}\hat{\hat{p}}_{\delta,\gamma} - p_\delta\right\|_\infty = \left\|\mathbb{E}\hat{\hat{p}}_{\delta,\gamma} - \mathbb{E}\hat{p}_\delta\right\|_\infty \leq \sup_x \mathbb{E}|\hat{\hat{p}}_{\delta,\gamma}(x) - \hat{p}_\delta(x)| \leq \frac{z'\varepsilon}{\delta}$ on $E_{m,\varepsilon}$.

Let $\varepsilon_m = c''\sqrt{\frac{\tau - \log n}{m\gamma_m}}$, then applying Theorem **5.3.1.** gives the desired result.

# References

Aghabozorgi, S., Seyed Shirkhorshidi, A., Ying Wah, T. (2015) Time-series clustering – A decade review. *Information systems (Oxford).* [Online] 5316–38.

Ashkartizabi, M. & Aminghafari, M. (2017) Functional data clustering using K-means and random projection with applications to climatological data. *Stochastic environmental research and risk assessment.* [Online] 32 (1), 83–104.

Ben-Hur, A. & Guyon, I. (2003) Detecting stable clusters using principal component analysis. *Methods in molecular biology* (Clifton, N.J.). 224159–182.

Berrendero, J., Justel, A., Svarc, M. (2011) Principal components for multivariate functional data. *Computational statistics & data analysis.* [Online] 55 (9), 2619–2634.

Chaudhuri, K. & Dasgupta, S. (2014) Consistent Procedures for Cluster Tree Estimation and Pruning. *IEEE transactions on information theory.* [Online] 60 (12), 7900–7912.

Cuevas, A. & Rodríguez-Casal, A. (2004) On boundary estimation. *Advances in applied probability.* [Online] 36 (2), 340–354.

Dasgupta, S. & Freund, Y. (2009) Random Projection Trees for Vector Quantization. *IEEE transactions on information theory.* [Online] 55 (7), 3229–3242.

Devroye, L. and Lugosi, G. (2001) Combinatorial Methods in Density Estimation. Springer, New York.

Ester, M., Kriegel, H.-P., Sander, J., Xu, X. (1996) A density-based algorithm for discovering clusters in large spatial databases with noise. In: Simoudis, E., Han, J., Fayyad, U. (eds): *Proceedings of the Second International Conference on Knowledge Discovery and Data Mining* (KDD-96), 226–231., AAAI Press

Fuchs, E., Gruber, T., Nitschke, J., Sick, B. (2009) On-line motif detection in time series with SwiftMotif. *Pattern recognition.* [Online] 42 (11), 3015–3031.

Giné, E. & Guillou, A. (2001) On consistency of kernel density estimators for randomly censored data: rates holding uniformly over adaptive intervals. *Annales de l'I.H.P. Probabilités et statistiques*. [Online] 37 (4), 503–522.

Hartigan, J. A. (1977) Distribution Problems in Clustering. In: *Classification and Clustering.* s.l.:Academic Press, pp. 45-71.

Hartigan, J. A. (1981) Statistical Theory in Clustering. *Journal of Classification,* pp. 63-76.

Jacques, J. & Preda, C. (2013) Functional data clustering: a survey. *Advances in data analysis and classification.* [Online] 8 (3), 231–255.

Jacques, J. & Preda, C. (2014) Model-based clustering for multivariate functional data. *Computational statistics & data analysis*. [Online] 7192–106.

Jones, M. & Rice, J. (1992) Displaying the important features of large collections of similar curves. *The American statistician*. [Online] 46 (2), 140–145.

Li, H., Deng, X., Dolloff, C. A., Smith, E. P. (2016) Bivariate functional data clustering: grouping streams based on a varying coefficient model of the stream water and air temperature relationship. *Environmetrics* (London, Ont.). [Online] 27 (1), 15–26.

Li, S. (2011) Concise Formulas for the Area and Volume of a Hyperspherical Cap. *Asian Journal of Mathematics & Statistics*. 4: 66-70.

Maier, M., Hein, M., von Luxburg, U. (2009) Optimal construction of k -nearest-neighbor graphs for identifying noisy clusters. *Theoretical computer science*. [Online] 410 (19), 1749–1764.

Malik, S., Kanwal, N., Asghar, M. N., et al. (2019) Data Driven Approach for Eye Disease Classification with Machine Learning. *Applied sciences*. [Online] 9 (14), 2789–.

Park, J. & Ahn, J. (2017) Clustering multivariate functional data with phase variation. *Biometrics*. [Online] 73 (1), 324–333.

Polonik, W. (1995) Measuring Mass Concentrations and Estimating Density Contour Clusters-An Excess Mass Approach. *The Annals of statistics*. [Online] 23 (3), 855–881.

Ramsay, J. O. & Dalzell, C. J. (1991) Some Tools for Functional Data Analysis. *Journal of the Royal Statistical Society*. Series B, Methodological. [Online] 53 (3), 539–572.

Rigollet, P. & Vert, R. (2009) Optimal rates for plug-in estimators of density level sets. *Bernoulli : official journal of the Bernoulli Society for Mathematical Statistics and Probability*. [Online] 15 (4), 1154–1178.

Rinaldo, A. & Wasserman, L. (2010) Generalized Density Clustering. *The Annals of statistics*. [Online] 38 (5), 2678–2722.

Schmutz, A., Jacques, J., Bouveyron, C., Chèze, L., Martin, P. (2020) Clustering multivariate functional data in group-specific functional subspaces. *Computational statistics*. [Online] 35 (3), 1101–1131.

Singh, A., Scott, C., Nowak, R. (2009) Adaptive Hausdorff Estimation of Density Level Sets. *The Annals of statistics*. [Online] 37 (5B), 2760–2782.

Sriperumbudur, B. & Steinwart, I. (2012) Consistency and Rates for Clustering with DBSCAN. *Proceedings of the Fifteenth International Conference on Artificial Intelligence and Statistics*, in *Proceedings of Machine Learning Research* 22:1090-1098 Available from https://proceedings.mlr.press/v22/sriperumbudur12.html

Steinwart, I. (2015) Fully Adaptive Density-Based Clustering. *The Annals of statistics.* [Online] 43 (5), 2132–2167.

Talagrand, M. (1995) Concentration of measure and isoperimetric inequalities in product spaces. *Publications mathématiques. Institut des hautes études scientifiques*. [Online] 81 (1), 73–205.

Talagrand, M. (1996) New concentration inequalities in product spaces. *Inventiones mathematicae*. [Online] 126 (3), 505–563.

Tsybakov, A. B. (1997) On Nonparametric Estimation of Density Level Sets. *The Annals of statistics*. [Online] 25 (3), 948–969.

Wang, J.-L., Chiou, J.-M., Mueller, H.-G. (2015). Review of Functional Data Analysis.

Zambom, A. Z., Collazos, J., Dias, R. (2018) Functional data clustering via hypothesis testing k-means. *Computational statistics*. [Online] 34 (2), 527–549.